\documentclass[11pt]{article}
\usepackage{amsfonts}
\date{\empty}

\usepackage{amsmath}
\usepackage{amssymb}
\usepackage{latexsym}
\newcommand \nc {\newcommand}
\nc \proof {\noindent {\em{Proof.\/ }}}
\nc \qed {\hfill $\Box$}

\newtheorem{Theorem}{Theorem}[section]
\newtheorem{Lemma}{Lemma}[section]

\newtheorem{definition}[Theorem]{Definition}

\newcommand{\beq}{\begin{eqnarray}}
\newcommand{\eeq}{\end{eqnarray}}
\newcommand{\beqno}{\begin{eqnarray*}}
\newcommand{\eeqno}{\end{eqnarray*}}
\newcommand{\be}{\begin{equation}}
\newcommand{\ee}{\end{equation}}

\def \d2{\Delta_{+}^2}

\begin{document}

\title{Well-posedness for the heat flow of biharmonic maps with rough initial data}
\author{Changyou Wang
\thanks{Department of Mathematics, University of Kentucky,
Lexington, KY 40506}}
\date{}
\maketitle
\begin{abstract}
This paper establishes the local (or global, resp.)  well-posedness of the heat flow
of bihharmonic maps from $\mathbb R^n$ to a compact Riemannian manifold without boundary
for initial data with small local BMO (or BMO, resp.)  norms.

\end{abstract}
\setcounter{section}{0} \setcounter{equation}{0}
\section{Introduction}

For $k\ge 1$, let $N$ be a $k$-dimensional
compact Riemannian manifold without boundary, isometrically
embedded in some Euclidean space $\mathbb R^l$.
Let $\Omega\subset\mathbb R^n$, $n\ge 1$, be a smooth domain.
There are two second order energy functional for mappings
from $\Omega$ to $N$, namely, the Hessian energy functional and
tension field energy functional given by
$$F(u)=\int_\Omega |\Delta u|^2, \  E(u)=\int_\Omega |D\Pi(u)(\Delta u)|^2,
\ u\in W^{2,2}(\Omega, N),$$
where $\Pi: N_{\delta_N}\to N$ is the smooth nearest point projection from
$N_{\delta_N}=\{y\in\mathbb R^l: \ {\rm{dist}}(y,N)\le\delta_N\}$ to $N$ 
for some small $\delta_N>0$, and 
$$W^{2,2}(\Omega,N)=\{v\in W^{2,2}(\Omega, \mathbb R^l): \ v(x)\in N
\ {\rm{for\ a.e.}}\ x\in\Omega\}.$$

Recall that a map $u\in W^{2,2}(\Omega,N)$ is called an (extrinsic) biharmonic map
(or intrinsic biharmonic map, resp.) if $u$ is a critical point of $F(\cdot)$
(or $E(\cdot)$, resp.).  Geometrically, a biharmonic map $u$ to $N$ enjoys the 
property that $\Delta^2 u$ is perpendicular to $T_u N$. The Euler-Lagrange equation
for biharmonic maps (see \cite{W2}) is:
\beq{} \label{biharm_map}
\Delta^2 u=\Delta(D^2\Pi(u)(\nabla u,\nabla u))
+2\nabla\cdot\langle\Delta u, \nabla(D\Pi(u))\rangle
-\langle \Delta u, \Delta(D\Pi(u))\rangle. 
\eeq
The Euler-Lagrange equation for intrinsic biharmonic maps (see \cite{W2}) is:
\begin{eqnarray} \label{int_biharm_map}
\Delta^2 u&=&\Delta(D^2\Pi(u)(\nabla u,\nabla u))
+2\nabla\cdot\langle\Delta u, \nabla(D\Pi(u))\rangle
-\langle \Delta u, \Delta(D\Pi(u))\rangle\nonumber\\
&+&D\Pi(u)[D^2\Pi(u)(\nabla u,\nabla u)\cdot D^3\Pi(u)(\nabla u,\nabla u)]\nonumber\\
&+& 2D^2\Pi(u)(\nabla u,\nabla u)\cdot D^2\Pi(u)(\nabla u, \nabla(D\Pi(u))).
\end{eqnarray}

The study of biharmonic maps was initiated by Chang-Wang-Yang \cite{CWY} in late 90's.
It has since drawn considerable research interests. In particular, the smoothness of 
biharmonic maps  (and intrinsic biharmonic maps)  in $W^{2,2}$ has been established in dimension $4$ by \cite{CWY} for $N=S^{l-1}$ and by \cite{W1} for general manifold $N$. For $n\ge 5$, the partial regularity of the class of stationary biharmonic maps in
$W^{2,2}$ has been shown by by \cite{CWY} for $N=S^{l-1}$ and by \cite{W1} for general manifold $N$.  The readers can refer to Strzelecki \cite{Str},  Angelesberg \cite{A},
Lamm-Riviere \cite{LR},
Struwe \cite{Struwe}, Scheven \cite{Sch}, Hong-Wang \cite{HW}, and Wang \cite{W3}
for further interesting results. 

Motivated by the study of heat flow of harmonic maps, which has played a very important role
in the existence of harmonic maps in various topological classes, it is very natural and interesting
to study the corresponding heat flow of biharmonic maps.    For $\Omega=\mathbb R^n$,
the heat flow of harmonic maps for $u:\mathbb R^n\times \mathbb R_+ \to N$
is given by 
\begin{eqnarray} \label{hfbhm}
\partial_t u+\Delta^2 u&=&\Delta(D^2\Pi(u)(\nabla u,\nabla u))+2\nabla\cdot\langle\Delta u, \nabla(D\Pi(u))\rangle\nonumber\\
&&-\langle\nabla\Delta u, \Delta(D\Pi(u))\rangle \ \ 
\ \ \ \ \ \ \ \ \ \ \ \ \ \  \  \ \ \ {\rm{in}}\ \mathbb R^n\times (0,+\infty)\\
u|_{t=0}&=& u_0 \ \ \ \ \ \ \ \ \ \ \ \ \  \ \ \ \ \ \ \ \ \ \ \ \ \ \ \ \ \ \  \ \  \ \ \ \ \  \ \ \ \ \ \ \ {\rm{on}}\ \mathbb R^n,   \label{ic}
\end{eqnarray}
where $u_0:\mathbb R^n\to N$ is a given map.

(\ref{hfbhm})-(\ref{ic}) was first investigated by Lamm in \cite{L1, L2}, 
where for smooth initial data $u_0\in C^\infty(\mathbb R^n,N)$  the short time
smooth solution was established. 
Moreover, such  a short time smooth solution is proven to be  globally smooth provided that
$n=4$ and $\|u_0\|_{W^{2,2}(\mathbb R^4)}$ is sufficiently small.  For large initial data
$u_0\in W^{2,2}(\mathbb R^4)$, it was independently proved by Gastel \cite{G} and
Wang \cite{W4} that there exists a global weak solution to (\ref{hfbhm})-(\ref{ic})
that is smooth away from finitely many singular times. 

It is a very interesting question to  seek the largest class of rough initial data such that
(\ref{hfbhm})-(\ref{ic}) is well-posed (either local or global) in suitable spaces.  There have
been interesting works on this type of question for the Navier-Stokes equation (see
Koch-Tataru \cite{KT}), the heat flow of harmonic maps (see Koch-Lamm \cite{KL} and
Wang \cite{W5}), and the Willmore flow, the Ricci flow, and the Mean curvature flow by
Koch-Lamm \cite{KL}. 

The main goal of this paper is to investigate the well-posedness issue
of (\ref{hfbhm}) and (\ref{ic}) for initial data $u_0$ with small BMO norm.

To state our main result, we first introduce the BMO spaces. 
\begin{definition} For $0<R\le+\infty$,
the local BMO space, ${\rm{BMO}}_R(\mathbb R^n)$, is the space consisting of
locally integrable functions $f$ such that 
$$
\left[f\right]_{\rm{BMO}_R(\mathbb R^n)}
:=\sup_{x\in\mathbb R^n, 0<r\le R} \{r^{-n}\int_{B_r(x)}|f-f_{x,r}| \}<+\infty,
$$
where  $B_r(x)\subset\mathbb R^n$ is the ball  with center $x$ and radius $r$,
and  
$$f_{x,r}=\frac{1}{|B_r(x)|}\int_{B_r(x)} f $$ is the average of $f$ over $B_r(x)$. 
We say $f\in {\overline{\rm VMO}}(\mathbb R^n)$ if
$$\lim_{r\downarrow 0}\left[f\right]_{{\rm BMO}_r(\mathbb R^n)}=0.$$
For $R=+\infty$, we simply write 
$({\rm{BMO}}(\mathbb R^n), [\cdot]_{{\rm{BMO}}(\mathbb R^n)})$
for $({\rm{BMO}}_\infty(\mathbb R^n), [\cdot]_{{\rm{BMO}}_\infty(\mathbb R^n)})$.
\end{definition}

For $0<T\le +\infty$, we also introduce the functional space $X_{T}$ as follows.
\begin{equation}\label{space1}
X_T=\left\{f:\mathbb R^n\times [0,T]\to \mathbb R\ | \
\|f\|_{X_T}\equiv \sup_{0<t\le T}\|f(t)\|_{L^\infty(\mathbb R^n)}+\left[f\right]_{X_T}<+\infty\right\}
\end{equation}
where
\begin{eqnarray}\label{x}
\left[f\right]_{X_T}&=&\sup_{0<t\le T}(\sum_{i=1}^2 t^{\frac{i}4}\|\nabla^i f(t)\|_{L^\infty(\mathbb R^n)})
+\sup_{x\in\mathbb R^n, 0<R\le T^\frac14} \ (R^{-n}\int_{P_R(x,R^4)}|\nabla f|^4)^\frac14\nonumber\\
&&+\sup_{x\in\mathbb R^n, 0<R\le T^\frac14} \ (R^{-n}\int_{P_R(x,R^4)}|\nabla^2 f|^2)^\frac12,
\end{eqnarray}
where $P_R(x,R^4)=B_R(x)\times [0,R^4]$ is the parabolic cylinder with
center $(x,R^4)$ and radius $R$. 
It is clear that $(X_T, \|\cdot\|_{X_T})$ is a Banach space.
When $T=+\infty$, we simply write $X$ for $X_\infty$, $\|\cdot\|_X$ for
$\|\cdot\|_{X_\infty}$, and $[\cdot]_X$ for $[\cdot]_{X_\infty}$.

The first theorem states
\begin{Theorem}\label{gwp} There exists $\epsilon_0>0$ such that
for any $R>0$ if $[u_0]_{{\rm BMO}_R(\mathbb R^n)}\le\epsilon_0$, then there exists 
a unique solution $u\in X_{R^4}$ to
(\ref{hfbhm})-(\ref{ic}) with small $[u]_{X_{T}}$. In particular, if $u_0
\in{\overline{\rm VMO}}(\mathbb R^n)$ then there exists $T_0>0$ such that
(\ref{hfbhm})-(\ref{ic}) admits a unique solution $u\in X_{T_0}$ with small
$[u]_{X_{T_0}}$.
\end{Theorem}

As a direct corollary, we have the following global well-posedness result.
\begin{Theorem} There exists $\epsilon_0>0$ such that
if $[u_0]_{{\rm BMO}(\mathbb R^n)}\le\epsilon_0$, then there exists 
a unique solution $u\in X$ to
(\ref{hfbhm})-(\ref{ic}) with small $[u]_{X}$.
\end{Theorem}

Now we turn to the discussion of the heat flow of intrinsic biharmonic maps.
The equation of the heat flow of  intrinsic biharmonic maps on $\mathbb R^n$
is given by
\begin{eqnarray} \label{int_bih}
\partial_t u+\Delta^2 u&=&\Delta(D^2\Pi(u)(\nabla u,\nabla u))
+2\nabla\cdot\langle\Delta u, \nabla(D\Pi(u))\rangle
-\langle \Delta u, \Delta(D\Pi(u))\rangle\nonumber\\
&+&D\Pi(u)[D^2\Pi(u)(\nabla u,\nabla u)\cdot D^3\Pi(u)(\nabla u,\nabla u)]\nonumber\\
&+& 2D^2\Pi(u)(\nabla u,\nabla u)\cdot D^2\Pi(u)(\nabla u, \nabla(D\Pi(u))
\  {\rm in}\ \mathbb R^n\times (0,+\infty)\\
u\big|_{t=0}&=& u_0: \mathbb R^n\to N. \label{ic1}
\end{eqnarray}
In \cite{L3}, Lamm studied (\ref{int_bih})-(\ref{ic1}).  Under the assumption that $n\le 4$
and the section curvature of $N$ is nonpositive, the global smooth solution
to (\ref{int_bih})-(\ref{ic1}) was established in \cite{L3}.  

Analogous to Theorem 1.2 and 1.3, we obtain the following results on (\ref{int_bih})-(\ref{ic1}).
\begin{Theorem} There exists $\epsilon_0>0$ such that
for any $R>0$ if $[u_0]_{{\rm BMO}_R(\mathbb R^n)}\le\epsilon_0$, then there exists 
a unique solution $u\in X_{R^4}$ to
(\ref{int_bih})-(\ref{ic1}) with small $[u]_{X_{T}}$. In particular, if $u_0
\in{\overline{\rm VMO}}(\mathbb R^n)$ then there exists $T_0>0$ such that
(\ref{int_bih})-(\ref{ic1}) admits a unique solution $u\in X_{T_0}$ with small
$[u]_{X_{T_0}}$.
\end{Theorem}

\begin{Theorem} There exists $\epsilon_0>0$ such that
 if $[u_0]_{{\rm BMO}(\mathbb R^n)}\le\epsilon_0$, then there exists 
a unique solution $u\in X$ to
(\ref{int_bih})-(\ref{ic1}) with small $[u]_{X}$.
\end{Theorem}

We remark that since $W^{1,n}(\mathbb R^n)\subset {\overline{\rm VMO}}(\mathbb R^n)$,
it follows from Theorem 1.2 (or Theorem 1.4, resp.) that (\ref{hfbhm})-(\ref{ic}) 
(or {\ref{int_bih})-(\ref{ic1}), resp.) is uniquely solvable
in $X_{T_0}$ for some $T_0>0$  provided $u_0\in W^{1,n}(\mathbb R^n, N)$;
and is uniquely solvable in $X$ provided $\|\nabla u_0\|_{L^{n}(\mathbb R^n)}$ is sufficiently small, via Theorem 1.3 (or Theorem 1.5, resp.).

We also remark that the techniques to handle the heat flow of biharmonic maps illustrated in this paper can be extended
to investigate the well-posedness of the heat flow of polyharmonic maps for BMO initial data in any dimensions.
This will be discussed in a forthcoming paper \cite{HW1}.

The remaining of the paper is written as follows. In section 2, we review some basic estimates
on the biharmonic heat kernel, due to Koch-Lamm \cite{KL}.  In section 3, we outline some
crucial estimates on the biharmonic heat equation. In section 4, we prove the boundedness
of the mapping operator $\mathbb S$ determined by the Duhamel formula. In section 5,
we prove Theorem 1.2 and 1.3. In section 6, we prove Theorem 1.4 and 1.5. 

\setcounter{section}{1} \setcounter{equation}{0}
\section{Review of the biharmonic heat kernel}

In this section, we review some fundamental properties from Koch and Lamm \cite{KL} on the
biharmonic heat kernel.

Consider the fundamental solution of the biharmonic heat equation:
$$(\partial_t+\Delta^2)b(x,t)=0  \ {\rm{in}}\ \mathbb R^n\times \mathbb R_+$$
and it is given by
$$b(x,t)=t^{-\frac{n}4}g(\frac{x}{t^\frac14}), $$
where
\begin{equation}
\label{bikernel}
g(\xi)=(2\pi)^{-\frac{n}2}\int_{\mathbb R^n}e^{i\xi k-|k|^4}\,dk, \ \xi\in\mathbb R^n.
\end{equation}

The following Lemma, due to Koch and Lamm \cite{KL} (Lemma 2.4), play a very important role in this paper.
\begin{Lemma} For $x\in\mathbb R^n$ and $t>0$, the following estimates hold:
\begin{equation}\label{est1}
|b(x,t)|\le ct^{-\frac{n}4}\exp(-\alpha \frac{|x|^\frac43}{t^\frac13}), \ \alpha=\frac{3 2^\frac13}{16},
\end{equation}
\begin{equation}\label{est2}
|\nabla^k b(x,t)|\le c(t^\frac14+|x|)^{-n-k}, \ \forall k\ge 1
\end{equation}
\begin{equation}\label{est3}
\|\nabla^k b(\cdot,t)\|_{L^1(\mathbb R^n)}\le ct^{-\frac{k}{4}}, \ \forall k\ge 1.
\end{equation}
Moreover, there exist $c,c_1>0$ such that for $0\le j\le 4$,
\begin{equation}\label{est4}
|\nabla^j b(x,t)|\le ce^{-c_1|x|}, \ \forall (x,t)\in\mathbb R^n\times (0,1)\setminus (B_2\times (0,\frac12)).
\end{equation}
\end{Lemma}

For the purpose of this paper, we also recall the Carleson's characterization of BMO spaces.
Let $\mathcal S$ denote the class of Schwartz functions. Then the following property is
well-known (see, Stein \cite{Stein}).
\begin{Lemma} Let $\Phi\in\mathcal S$ be such that $\int_{\mathbb R^n}\Phi=0$. For $t>0$,
let $\Phi_t(x)=t^{-n}\Phi(\frac{x}{t}), \ x\in\mathbb R^n$. If $f\in\rm{BMO}(\mathbb R^n)$,
then $|\Phi_t*f|^2(x,t)\frac{dxdt}{t}$ is a Carleson measure on $\mathbb R^{n+1}_+$, i.e.,
\begin{equation}
\sup_{x\in\mathbb R^n, r>0} r^{-n}\int_0^r\int_{B_r(x)}|\Phi_t*f|^2\frac{dxdt}{t}
\leq C[u_0]_{\rm{BMO}(\mathbb R^n)}^2 \label{carleson}
\end{equation}
for some $C=C(n)>0$.
If $f\in {\rm BMO}_R(\mathbb R^n)$ for some $R>0$, then 
\begin{equation}
\sup_{x\in\mathbb R^n,0< r\le R} r^{-n}\int_0^r\int_{B_r(x)}|\Phi_t*f|^2\frac{dxdt}{t}
\leq C[u_0]_{{\rm BMO}_R(\mathbb R^n)}^2 \label{carleson1}
\end{equation}
for some $C=C(n)>0$.
\end{Lemma}

Recall that the solution to the Dirichlet problem of the inhomogeneous biharmonic heat equation
\begin{eqnarray}\label{biheat}
(\partial_t +\Delta^2)u&=& f \ \ \rm{on} \ \ \mathbb R^n\times(0,+\infty)\\
u&=& u_0 \ \ {\rm on} \ \ \mathbb R^n\times\{0\}
\end{eqnarray}
is given by the Duhamel formula:
\begin{equation}\label{duhemel}
u=\mathbb Gu_0+ \mathbb Sf
\end{equation}
where
\begin{equation}
\mathbb Gu_0(x,t):=(b(\cdot,t)*u_0)(x)=\int_{\mathbb R^n} b(x-y,t)u_0(y)\,dy,
\ (x,t)\in\mathbb R^n\times(0,+\infty), \label{duhemel1}
\end{equation}
and
\begin{equation}
\mathbb Sf(x,t)=\int_0^t\int_{\mathbb R^n} b(x-y,t-s)f(y,s)\,dyds, \
(x,t)\in\mathbb R^n\times(0,+\infty). \label{duhemel2}
\end{equation}

\setcounter{section}{2} \setcounter{equation}{0}
\section{Basic estimates for the biharmonic heat equation}

In this section, we provide some crucial estimates for the solution of the biharmonic heat equation with initial data in BMO spaces, including the estimate of the distance to the
manifold $N$.

\begin{Lemma} For $0<R\le +\infty$, if $u_0\in {\rm BMO}_R(\mathbb R^n)$, then $\hat{u}_0\equiv \mathbb Gu_0$ satisfies the
following estimates:
\begin{equation}\label{biheat_est1}
\sup_{x\in\mathbb R^n, 0<r\le R} r^{-n}\int_{P_r(x,r^4)}(|\nabla^2 \hat{u}_0|^2
+r^{-2}|\nabla \hat{u}_0|^2)\le C\left[u_0 \right]_{{\rm BMO}_R(\mathbb R^n)}^2,
\end{equation}
and
\begin{equation}\label{biheat_est2}
\sup_{0<t\le R^4}\left(\sum_{i=1}^2 t^{\frac{i}4}\|\nabla\hat{u}_0(t)\|_{L^\infty(\mathbb R^n)}\right)
\le C\left[u_0 \right]_{{\rm BMO}_R(\mathbb R^n)}.
\end{equation}
If, in addition, $u_0\in L^\infty(\mathbb R^n)$, then
\begin{equation}\label{biheat_est3}
\sup_{x\in\mathbb R^n, 0<r\le R} r^{-n}\int_{P_r(x,r^4)}|\nabla \hat{u}_0|^4
\le C\|u_0\|_{L^\infty(\mathbb R^n)}^2\cdot\left[u_0 \right]_{{\rm BMO}_R(\mathbb R^n)}^2.
\end{equation}
\end{Lemma}
\noindent{\it Proof}.  For simplicity, we present the argument for $R=+\infty$.
Let $g$ be given by (\ref{bikernel}).
Let $\Phi^i=\nabla^i g$ for $i=1,2$. Then it is clear that $\Phi^i\in\mathcal S$ and
$\int_{\mathbb R^n}\Phi^i=0$ for $i=1,2$. Hence by Lemma 2.2,
$|\Phi^i_t*u_0|^2\frac{dxdt}{t}$ is a Carleson measure on $\mathbb R^{n+1}_+$ for
$i=1,2$. Direct calculations show, for $i=1,2$,
$$\Phi^i_t(x)=t^{-n}(\nabla^i g)(\frac{x}{t})
=t^{i}\nabla^i\left(t^{-n}g(\frac{x}{t})\right)=t^i \nabla^i (g_t(x)),$$
where
$$g_t(x)=t^{-n}g(\frac{x}{t}).$$
Hence we have
$$(\Phi^i_t*u_0)(x)=t^i\nabla^i(g_t*u_0)(x).$$
Since the biharmonic heat kernel  $b(x,t)=g_{t^\frac14}(x)$, we have
$$(\Phi^i_t*u_0)(x)=t^i \nabla^i\left((b(\cdot, t^4)*u_0)(x)\right)
=t^i \nabla^i(\mathbb Gu_0)(x,t^4).$$
Thus we have, for $i=1,2$,
\begin{eqnarray*} C[u_0]_{\rm{BMO}(\mathbb R^n)}^2
&\geq&
\sup_{x\in\mathbb R^n, r>0} r^{-n}\int_0^r\int_{B_r(x)} |\Phi_t^i*u_0|^2\frac{dxdt}{t}\\
&=&\sup_{x\in\mathbb R^n, r>0} r^{-n}\int_0^r\int_{B_r(x)} t^{2i-1}|\nabla^i \mathbb Gu_0|^2(x,t^4)\,dxdt\\
&=&\frac14 \sup_{x\in\mathbb R^n, r>0} r^{-n}\int_{P_r(x,r^4)} t^{\frac{2i-4}4}
|\nabla^i \mathbb Gu_0|^2(x,t)\,dxdt
\end{eqnarray*}
This clearly implies (\ref{biheat_est1}), since for $i=1,2$,
$t^{\frac{2i-4}{4}}\ge r^{2i-4}$ when $0\le t\le r^4$.

Since $\hat{u}_0$ solves the biharmonic heat equation
$({\partial_t}+\Delta^2)\hat{u}_0=0$ on $\mathbb R^n\times(0,+\infty)$, the standard gradient estimate implies that for any $x\in\mathbb R^n$ and $r>0$,
$$
r^2|\nabla \hat{u}_0|^2(x,r^4)+r^4|\nabla^2 \hat{u}_0|^2(x,r^4)
\leq Cr^{-n}\int_{P_r(x,r^4)}(r^{-2}|\nabla \hat{u}_0|^2 +|\nabla^2 \hat{u}_0|^2).$$
Taking supremum over $x\in\mathbb R^n$ and setting $t=r^4>0$ yields (\ref{biheat_est2}).

For (\ref{biheat_est3}), observe that $u_0\in L^\infty(\mathbb R^n)$ implies
$\Phi^1_t*u_0\in L^\infty(\mathbb R^n)$ and
$$\|\Phi^1_t*u_0\|_{L^\infty(\mathbb R^n)}
\le \|\Phi^1\|_{L^1(\mathbb R^n)}\|u_0\|_{L^\infty(\mathbb R^n)}
\le \|\nabla g\|_{L^1(\mathbb R^n)}\|u_0\|_{L^\infty(\mathbb R^n)}\le C\|u_0\|_{L^\infty(\mathbb R^n)}.$$
Hence
\begin{eqnarray*}
&&\sup_{x\in\mathbb R^n, r>0}
\int_{P_r(x,r^4)} |\nabla \mathbb Gu_0|^4\,dxdt\\
&=&\sup_{x\in\mathbb R^n, r>0}
\int_0^r\int_{B_r(x)} |\Phi^1_t*u_0|^4\frac{dxdt}{t}\\
&\leq& \left(\sup_{t>0}\|\Phi^1_t*u_0\|_{L^\infty(\mathbb R^n)}\right)\cdot
\sup_{x\in\mathbb R^n, r>0}
\int_0^r\int_{B_r(x)} |\Phi^1_t*u_0|^2\frac{dxdt}{t}\\
&\leq& C\|u_0\|_{L^\infty(\mathbb R^n)}^2\cdot\left[u_0\right]_{{\rm{BMO}}(\mathbb R^n)}^2.
\end{eqnarray*}
This implies (\ref{biheat_est3}). \qed

Now we prove an important estimate on the distance of $\hat{u}_0$ to the manifold
$N$ in terms of the BMO norm of $u_0$. More precisely, 
\begin{Lemma} \label{distance} For any $\delta>0$, there exists $K=K(\delta,N)>0$ such that
for $R>0$ if $u_0\in{\rm BMO}_R(\mathbb R^n)$ then 
\beq{}\label{distance1}
{\rm{dist}}(\hat{u}_0(x,t),N)\le K \left[u_0\right]_{{\rm BMO}_R(\mathbb R^n)}+\delta,
\ \ \forall x\in\mathbb R^n, \ 0\le t\le \frac{R^4}{K^4}.
\eeq
In particular, if  $u_0\in{\rm BMO}(\mathbb R^n)$ then 
\beq{}\label{distance20}
{\rm{dist}}(\hat{u}_0(x,t),N)\le K \left[u_0\right]_{{\rm BMO}(\mathbb R^n)}+\delta,
\ \ \forall x\in\mathbb R^n, \ t\in\mathbb R_+.
\eeq
\end{Lemma}
\noindent{\it Proof}. Since (\ref{distance20}) follows directly from (\ref{distance1}),
it suffices to prove (\ref{distance1}). 
For any $x\in\mathbb R^n$, $t>0$, and $K>0$, denote
$$c_{x,t}^K=\frac{1}{|B_K(0)|}\int_{B_K(0)} u_0(x-t^\frac14 z)\,dz.$$
Let $g$ be given by (2.1). Then, by a change of variables, we have
$$\hat{u}_0(x,t)=\int_{\mathbb R^n} g(y) u_0(x-t^\frac14 y)\,dy.$$
Applying Lemma 2.1, we have
\begin{eqnarray}\label{distance2}
\left|\hat{u}_0(x,t)-c_{x,t}^K\right|
&\leq& \int_{\mathbb R^n} g(y) |u_0(x-t^\frac14 y)-c_{x,t}^K|\,dy
\nonumber\\
&\leq&\left\{\int_{B_K(0)}+\int_{\mathbb R^n\setminus B_K(0)} \right\}
g(y)|u_0(x-t^\frac14 y)-c_{x,t}^K|\,dy\nonumber\\
&\leq& \int_{B_K(0)}ce^{-\alpha |y|^\frac43}|u_0(x-t^\frac14 y)-c_{x,t}^K|\,dy\nonumber\\
&&+2\|u_0\|_{L^\infty(\mathbb R^n)}\int_{\mathbb R^n\setminus B_K(0)}c e^{-\alpha |y|^\frac43}\,dy\nonumber\\
&\leq& K^n \left[u_0\right]_{{\rm BMO}_{Kt^{\frac14}}(\mathbb R^n)}+C_N\int_K^\infty e^{-\alpha r^\frac43} r^{n-1}\,dr\nonumber\\
&\leq& \delta+ K^n \left[u_0\right]_{{\rm BMO}_{Kt^{\frac14}}(\mathbb R^n)}
\end{eqnarray}
provide we choose a sufficiently large $K=K(\delta,N)>0$ so that 
$$C_N\int_K^\infty e^{-\alpha r^\frac43} r^{n-1}\,dr\le\delta.$$
On the other hand, since $u_0(\mathbb R^n)\subset N$, we have
$$
{\rm{dist}}(c_{x,t}^K, N)\le \left|c_{x,t}^K-u_0(x-t^\frac14 y)\right|, \ \forall y\in B_K(0)$$
and hence
\begin{equation}\label{distance3}
{\rm{dist}}(c_{x,t}^K, N)\le \frac{1}{|B_K(0)|}\int_{B_K(0)}|c_{x,t}^K-u_0(x-t^\frac14 y)|\,dy
\le \left[u_0\right]_{{\rm BMO}_{Kt^\frac14}(\mathbb R^n)}.
\end{equation}
Putting (\ref{distance2}) and (\ref{distance3}) together  yields (\ref{distance1})
holds for $t\le \frac{R^4}{K^4}$.
This completes the proof. \qed

\setcounter{section}{3} \setcounter{equation}{0}
\section{Boundedness of the operator $\mathbb S$}

In this section, we introduce two more functional spaces and establish the boundedness of the operator $\mathbb S$ between these spaces.

For $0<T\le +\infty$, besides the space $X_T$ introduced in the section 1,
we need to introduce the spaces $Y_T^1, Y_T^2$.

The space $Y_T^1$ is the space consisting of functions $f:\mathbb R^n\times [0,T]\to\mathbb R$
such that
\begin{equation}\label{y1}
\|f\|_{Y_T^1}\equiv \sup_{0<t\le T}\ t\|f(t)\|_{L^\infty(\mathbb R^n)}
+\sup_{x\in\mathbb R^n, 0<r\le T^\frac14} \ r^{-n}\int_{P_r(x,r^4)}|f|<+\infty,
\end{equation}
and the space $Y_T^2$ is the space consisting of  functions $f:\mathbb R^n\times [0,T]\to\mathbb R$
such that
\begin{equation}\label{y2}
\|f\|_{Y_T^2}\equiv \sup_{0<t\le T}\ t^\frac34\|f(t)\|_{L^\infty(\mathbb R^n)}
+\sup_{x\in\mathbb R^n, 0<r\le T^\frac14} \ (r^{-n}\int_{P_r(x,r^4)}|f|^\frac43)^\frac34<+\infty.
\end{equation}
It is  easy to see $(Y_T^i, \|\cdot\|_{Y_T^i})$ is a Banach space for $i=1,2$.
When $T=+\infty$, we simply denote $(Y^i, \|\cdot\|_{Y^i})$ for $(Y_\infty^i, \|\cdot\|_{Y_\infty^i})$ 
for $i=1,2$.

Let the operator $\mathbb S$ be defined by (2.12). Then we have
\begin{Lemma}For any $0<T\le +\infty$,  if $f\in Y^1_T$, then $\mathbb Sf\in X_T$ and
\begin{equation}\label{bound1}
\|\mathbb Sf\|_{X_T}\le C\|f\|_{Y^1_T}
\end{equation}
for some $C=C(n)>0$.
\end{Lemma}
\noindent{\it Proof}.  We need to show the pointwise estimate
\begin{equation}\label{pt_est}
\sum_{i=0}^2 R^{i}|\nabla^i (\mathbb Sf)|(x,R^4)\le C\|f\|_{Y^1_T}, 
\ \forall x\in\mathbb R^n, 0<R\le T^\frac14,
\end{equation}
and the integral estimate for $0<R\le T^\frac14$:
\begin{equation}\label{int_est}
R^{-\frac{n}4}\|\nabla (\mathbb Sf)\|_{L^4(P_R(x,R^4))}
+R^{-\frac{n}2}\|\nabla^2 (\mathbb Sf)\|_{L^2(P_R(x,R^4))}\le C\|f\|_{Y^1_T} .
\end{equation}
By suitable scalings, we may assume $T\ge 1$.
Since both estimates are translation and scale invariant, it suffices to show that both
(\ref{pt_est}) and (\ref{int_est}) hold for $x=0$ and $R=1$.

For $i=0,1,2$, we have
\begin{eqnarray*}
\left|\nabla^i \mathbb Sf(0,1)\right|&=&\left|\int_0^1\int_{\mathbb R^n}\nabla^i b(y,1-s) f(y,s)\,dyds\right|\\
&\leq &\left\{\int_{\frac12}^1\int_{\mathbb R^n}+\int_0^{\frac12}\int_{B_2} +
\int_0^{\frac12}\int_{\mathbb R^n\setminus B_2}\right\} |\nabla^i b(y,1-s)| |f(y,s)|\,dyds\\
&=& I_1+I_2+I_3.
\end{eqnarray*}
Applying Lemma 2.1, we can estimate $I_1,I_2,I_3$ as follows.
\begin{eqnarray*}|I_1|&\leq& \left(\sup_{\frac12\le s\le 1}\|f(s)\|_{L^\infty(\mathbb R^n)}\right)
\left(\int_{\frac12}^1\|\nabla^i b(\cdot,1-s)\|_{L^1(\mathbb R^n)}\,ds\right)\\
&\leq& C\|f\|_{Y^1_1}\int_0^\frac12 s^{-\frac{i}4}\,ds\\
&\leq& C\|f\|_{Y^1_1}.
\end{eqnarray*}
\begin{eqnarray*}
|I_2|&\leq& (\sup_{0\le s\le\frac12}\|\nabla^i b(\cdot,1-s)\|_{L^\infty(\mathbb R^n)})
(\int_{B_2\times [0,\frac12]}|f(y,s)|\,dyds)\\
&\leq& C\int_{B_2\times [0,\frac12]}|f(y,s)|\,dyds\leq C\|f\|_{Y^1_1},
\end{eqnarray*}
and
\begin{eqnarray*}
|I_3|&\leq&\int_0^\frac12 \int_{\mathbb R^n\setminus B_2}|\nabla^i b(y,1-s)||f(y,s)|\,dyds\\
&\leq& C\int_0^\frac12\int_{\mathbb R^n\setminus B_2} e^{-c_1|y|}|f(y,s)|\,dyds\\
&\leq& C \left(\sum_{k=2}^\infty k^{n-1} e^{-c_1k}\right)
\cdot \left(\sup_{y\in \mathbb R^n} \int_{P_1(y,1)}|f(y,s)|\,dyds\right)\\
&\leq& C \|f\|_{Y^1_1}.
\end{eqnarray*}

Now we want to show (\ref{int_est}) by the energy method. Denote $w=\mathbb Sf$.
Then $w$ solves
\begin{equation} \label{w_eqn}
({\partial_t}+\Delta^2)w=f \ \ {\rm{in}}\  \mathbb R^n\times (0,+\infty);
\ w|_{t=0}=0.
\end{equation}
Let $\eta\in C^\infty_0(B_2)$ be a cut-off function of $B_1$. Multiplying (\ref{w_eqn}) by
$\eta^4 w$ and integrating over $\mathbb R^n\times [0,1]$, we obtain
$$\int_{\mathbb R^n\times \{1\}}|w|^2\eta^4
+2\int_{\mathbb R^n\times [0,1]}\Delta w\cdot\Delta(w\eta^4)=\int_{\mathbb R^n\times [0,1]}f\cdot w\eta^4.$$
This easily implies
\begin{eqnarray}\label{int_est10}
&&\int_{P_1(0,1)}|\nabla^2 w|^2\nonumber\\
&\leq&\int_{\mathbb R^n\times [0,1]}|\nabla^2 (\eta^2 w)|^2\nonumber\\
&\leq& C
\int_{\mathbb R^n\times [0,1]}\left[|\nabla\eta|^2|\nabla w|^2+(|\Delta\eta|+|\nabla\eta|^2)|w|^2\right]
+ C\int_{\mathbb R^n\times [0,1]} |f||w|\eta^2\nonumber\\
&\leq& C[\int_{(B_2\setminus B_1)\times [0,1]} |\nabla w|^2+|w|^2+\|f\|_{L^1(B_2\times [0,1])}
\|w\|_{L^\infty(B_2\times [0,1])}]\nonumber\\
&\leq & C \left[(\int_0^1 t^\frac12\,dt)\cdot (\sup_{0<t\le 1} t^{\frac12}\|\nabla w(t)\|_{L^\infty(\mathbb R^n)}^2)
+\|w\|_{L^\infty(B_2\times [0,1])}^2+\|f\|_{L^1(B_2\times [0,1])}^2\right]\nonumber\\
&\leq& C\left[\sup_{0<t\le 1} (\|w(t)\|_{L^\infty(\mathbb R^n)}^2+ t^\frac12\|\nabla w(t)\|_{L^\infty(\mathbb R^n)}^2)
+\|f\|_{Y^1_1}^2\right]\nonumber\\
&\leq& C\|f\|_{Y^1_1}^2,
\end{eqnarray}
where we have used (\ref{pt_est}) in the last step.

For the $L^4$ norm of $\nabla w$ on $P_1(0,1)$, recall the Nirenberg inequality  implies
$$\|\nabla(\eta^2 w(t))\|_{L^4(\mathbb R^n)}^4\leq C\|\eta^2 w(t)\|_{L^\infty(\mathbb R^n)}^2
\|\nabla^2(\eta^2 w(t))\|_{L^2(\mathbb R^n)}^2.$$
Integrating with respect to $t\in [0,1]$ clearly implies
$$\left(\int_{P_1(0,1)}|\nabla w|^4\right)^\frac14 \leq C \sup_{0\le t\le 1}\|w(t)\|_{L^\infty(\mathbb R^n)}^\frac12
\|\nabla^2(\eta^2 w)\|_{L^2(\mathbb R^n\times [0,1])}^{\frac12}\le
C\|f\|_{Y^1_1},$$
where we have used both (\ref{pt_est}) and (\ref{int_est10}) in the last step. This completes the proof.
\qed

To handle the nonlinearities of the heat flow of biharmonic maps (\ref{hfbhm}), we also need 

\begin{Lemma} For $0<T\le +\infty$,  if $f\in Y^2_T$, then for any $1\le\alpha\le n$,
$\mathbb S(\frac {\partial f}{\partial x_\alpha})\in X_T$ and
\begin{equation}\label{bound2}
\left\|\mathbb S(\frac {\partial f}{\partial x_\alpha})\right\|_{X_T}
\le C\left\|f\right\|_{Y^2_T}
\end{equation}
for some $C=C(n)>0$.
\end{Lemma}

\noindent{\it Proof}. 
The proof of (\ref{bound2}) is similar to that of Lemma 4.1. 
We will prove that for any $x\in\mathbb R^n$ and $0<R\le T^\frac14$,
both the pointwise estimate:
\begin{equation}\label{pt_est1}
\sum_{i=0}^2 R^{i}\left|\nabla^i (\mathbb S(\frac{\partial f}{\partial x_\alpha}))\right|(x,R^4)\le C\|f\|_{Y^2_T},
\end{equation}
and the integral estimate:
\begin{equation}\label{int_est1}
R^{-\frac{n}4}\left\|\nabla (\mathbb S(\frac{\partial f}{\partial x_\alpha}))\right\|_{L^4(P_R(x,R^4))}
+R^{-\frac{n}2}\left\|\nabla^2 (\mathbb S(\frac{\partial f}{\partial x_\alpha}))\right\|_{L^2(P_R(x,R^4))}\le C\|f\|_{Y^2_T}.
\end{equation}
By suitable scalings, we assume $T\ge 1$.
Since both estimates are translation and scale invariant, it suffices to show that both
(\ref{pt_est1}) and (4.10) hold for $x=0$ and $R=1$.
For $1\le\alpha\le n$, write $W_\alpha=\mathbb S(\frac{\partial f}{\partial x_\alpha})$. For $i=0,1,2$, we have
\begin{eqnarray*}
\nabla^i W_\alpha(0,1)&=&\int_{\mathbb R^n\times [0,1]} \nabla^i b(-y,1-s) \frac{\partial f}{\partial y_\alpha}(y,s)\,dyds\\
&=&\int_{\mathbb R^n\times [0,1]}(\nabla^i \frac{\partial}{\partial y_\alpha} b)(-y,1-s) f(y,s)\,dyds,
\end{eqnarray*}
which implies
\begin{eqnarray*}
\left|\nabla^i W_\alpha(0,1)\right|&\leq &\int_0^1\int_{\mathbb R^n}|\nabla^{i+1} b(y,1-s)| |f(y,s)|\,dyds\\
&=&\left\{\int_{\frac12}^1\int_{\mathbb R^n}+\int_0^{\frac12}\int_{B_2} +
\int_0^{\frac12}\int_{\mathbb R^n\setminus B_2}\right\} |\nabla^{i+1} b(y,1-s)| |f(y,s)|\,dyds\\
&=& I_4+I_5+I_6.
\end{eqnarray*}
Applying Lemma 2.1, we can estimate $I_4,I_5,I_6$ as follows.
\begin{eqnarray*}|I_4|&\leq& \left(\sup_{\frac12\le s\le 1}\|f(s)\|_{L^\infty(\mathbb R^n)}\right)
\left(\int_{\frac12}^1\|\nabla^{i+1} b(\cdot,1-s)\|_{L^1(\mathbb R^n)}\,ds\right)\\
&\leq& C\|f\|_{Y^2_1}\int_0^\frac12 s^{-\frac{i+1}4}\,ds\\
&\leq& C\|f\|_{Y^2_1},
\end{eqnarray*}
where we have used the fact  $\int_0^\frac12 s^{-\frac{i+1}4}\,ds<+\infty$ for $i\le 2$.
\begin{eqnarray*}
|I_5|&\leq& (\sup_{0\le s\le\frac12}\|\nabla^{i+1} b(\cdot,1-s)\|_{L^\infty(\mathbb R^n)})
(\int_{B_2\times [0,\frac12]}|f(y,s)|\,dyds)\\
&\leq& C\int_{B_2\times [0,\frac12]}|f(y,s)|\,dyds\leq C\|f\|_{Y^2_1},
\end{eqnarray*}
and since $i+1\le 3$, we have
\begin{eqnarray*}
|I_6|&\leq&\int_0^\frac12 \int_{\mathbb R^n\setminus B_2}|\nabla^{i+1} b(y,1-s)||f(y,s)|\,dyds\\
&\leq& C\int_0^\frac12\int_{\mathbb R^n\setminus B_2} e^{-c_1|y|}|f(y,s)|\,dyds\\
&\leq& C \left(\sum_{k=2}^\infty k^{n-1} e^{-c_1k}\right)
\cdot \left(\sup_{y\in \mathbb R^n} \int_{P_1(y,1)}|f(y,s)|\,dyds\right)\\
&\leq& C \|f\|_{Y^2_1}.
\end{eqnarray*}
Putting together these estimates, we prove (\ref{pt_est1}). (\ref{int_est1}) can be done by the energy method as well.
In fact, $W_\alpha$ solves
\begin{equation} \label{W_eqn}
({\partial_t}+\Delta^2)W_\alpha=\frac{\partial f}{\partial x_\alpha} \ \ {\rm{in}}\  \mathbb R^n\times (0,+\infty);
\ W_\alpha|_{t=0}=0.
\end{equation}
Let $\eta\in C^\infty_0(B_2)$ be a cut-off function of $B_1$. Multiplying (\ref{W_eqn}) by
$\eta^4 W_\alpha$ and integrating over $\mathbb R^n\times [0,1]$, we obtain
$$\int_{\mathbb R^n\times \{1\}}|W_\alpha|^2\eta^4
+2\int_{\mathbb R^n\times [0,1]}\Delta W_\alpha \cdot\Delta(W_\alpha\eta^4)
=-\int_{\mathbb R^n\times [0,1]}f\cdot \frac{\partial}{\partial x_\alpha}(W_\alpha\eta^4).$$
This implies
\begin{eqnarray}\label{int_est2}
&&\int_{P_1(0,1)}|\nabla^2 W_\alpha|^2\nonumber\\
&\leq&\int_{\mathbb R^n\times [0,1]}|\nabla^2 (\eta^2 W_\alpha)|^2\nonumber\\
&\leq& C
\int_{\mathbb R^n\times [0,1]}\left[|\nabla\eta|^2|\nabla W_\alpha|^2+(|\Delta\eta|+|\nabla\eta|^2)|W_\alpha|^2\right]\nonumber\\
&+& C\int_{\mathbb R^n\times [0,1]} |f|(|\nabla (\eta^2W_\alpha)|+|W_\alpha||\nabla \eta|)\nonumber\\
&\leq& C[\int_{(B_2\setminus B_1)\times [0,1]} (|\nabla W_\alpha|^2+|W_\alpha|^2)
+\|f\|_{L^1(B_2\times [0,1])}\|W_\alpha\|_{L^\infty(\mathbb R^n)}]\nonumber\\
&+& C\|f\|_{L^\frac43(B_2\times [0,1])}
\|\nabla (\eta^2 W_\alpha)\|_{L^4(\mathbb R^n\times [0,1])}\nonumber\\
&=&I_7+I_8.
\end{eqnarray}
It is easy to see that
\begin{eqnarray*}
&&|I_7|\\
&\leq& C\left[(\int_0^1 t^\frac12\,dt)\cdot (\sup_{0<t\le 1} t^{\frac12}\|\nabla W_\alpha(t)\|_{L^\infty(\mathbb R^n)}^2)
+\|W_\alpha\|_{L^\infty(B_2\times [0,1])}^2+\|f\|_{L^1(B_2\times [0,1])}^2\right]\\
&\leq& C\|f\|_{Y^2_1}^2
\end{eqnarray*}
where we have used the point wise estimate (\ref{pt_est1}) in the last step. In order to estimate $I_8$, we
 first need to employ the Nirenberg inequality:
$$\|\nabla(\eta^2 W_\alpha(t))\|_{L^4(\mathbb R^n)}^4\leq C\|\eta^2 W_\alpha(t)\|_{L^\infty(\mathbb R^n)}^2
\|\nabla^2(\eta^2 W_\alpha(t))\|_{L^2(\mathbb R^n)}^2,$$
which, after integrating with respect to $t\in [0,1]$, implies
\begin{equation}\label{int_est3}\|\nabla (\eta^2 W_\alpha)\|_{L^4(\mathbb R^n\times [0,1])}\leq C
\sup_{0\le t\le 1}\|W_\alpha(t)\|_{L^\infty(\mathbb R^n)}^\frac12
\|\nabla^2(\eta^2 W_\alpha)\|_{L^2(\mathbb R^n\times [0,1])}^{\frac12}.
\end{equation}
Therefore, $I_8$ can be estimated by
\begin{eqnarray*}
|I_8|&\leq& C \|f\|_{L^\frac43(B_2\times [0,1])}
\sup_{0\le t\le 1}\|W_\alpha(t)\|_{L^\infty(\mathbb R^n)}^\frac12
\|\nabla^2(\eta^2 W_\alpha)\|_{L^2(\mathbb R^n\times [0,1])}^{\frac12}\\
&\leq& \frac12 \int_{\mathbb R^n\times [0,1]}|\nabla^2(\eta^2 W_\alpha)|^2
+C\|f\|_{L^\frac43(B_2\times [0,1])}^\frac43
\sup_{0\le t\le 1}\|W_\alpha(t)\|_{L^\infty(\mathbb R^n)}^\frac23\\
&\leq& \frac12 \int_{\mathbb R^n\times [0,1]}|\nabla^2(\eta^2 W_\alpha)|^2
+C\|f\|_{L^\frac43(B_2\times [0,1])}^\frac43 \sup_{0\le t\le 1}\|W_\alpha(t)\|_{L^\infty(\mathbb R^n)}^\frac23\\
&\leq& \frac12 \int_{\mathbb R^n\times [0,1]}|\nabla^2(\eta^2 W_\alpha)|^2+C\|f\|_{Y^2_1}^2,
\end{eqnarray*}
where we have used (\ref{pt_est1}) in the last step. Now we substitute the estimates of $I_7$ and
$I_8$ into (\ref{int_est2}) and obtain
$$\int_{P_1(0,1)}|\nabla^2 W_\alpha|^2
\le C\int_{\mathbb R^n\times [0,1]}|\nabla^2(\eta^2 W_\alpha)|^2\leq C\|f\|_{Y^2_1}^2.$$
This, combined with (\ref{int_est3}), also implies
$$\int_{P_1(0,1)}|\nabla W_\alpha|^4 \le C\|f\|_{Y^2_1}^4.$$
The proof of (\ref{int_est1}) is now completed. \qed

\setcounter{section}{4} \setcounter{equation}{0}
\section{Proof of Theorem 1.2 and 1.3}

In this section, we will prove both Theorem 1.2 and 1.3. The idea is based on
the  fixed point theorem in a small ball inside $X_T$ for the mapping operator
determined by the Duhamel formula associate with (\ref{hfbhm})-(\ref{ic}).

First we need to extend $\Pi$ to $\mathbb R^l$. Let $\widetilde{\Pi}\in
C^\infty(\mathbb R^l,\mathbb R^l)$ be any smooth extension of $\Pi$ such that
$\widetilde{\Pi}\equiv\Pi$ on $N_{\delta_N}$.

Let
$$\mathcal F[u]= \Delta(D^2\widetilde\Pi(u)(\nabla u,\nabla u))
-\langle\Delta u, \Delta(D \widetilde\Pi(u))\rangle
+2\nabla\cdot\langle\Delta u, \nabla (D\widetilde\Pi(u)))\rangle
$$
be the right hand side nonlinearity of (\ref{hfbhm}).
Then it is easy to see that
\begin{eqnarray}
\mathcal F[u]&=&-\langle\Delta u, \Delta(D\widetilde\Pi(u))\rangle+
\nabla\cdot\left(2\langle\Delta u, \nabla(D\widetilde\Pi(u))\rangle
+\nabla(D^2\widetilde\Pi(u)(\nabla u,\nabla u))\right)\nonumber\\
&=& \mathcal F_1[u]+\nabla\cdot(\mathcal F_2[u]) \label{decom_non},
\end{eqnarray}
where
\begin{equation}\label{non}
\mathcal F_1[u]=-\langle\Delta u, \Delta(D\widetilde\Pi(u))\rangle,
\ \mathcal F_2[u]=2\langle\Delta u, \nabla(D\widetilde\Pi(u))\rangle
+\nabla(D^2\widetilde\Pi(u)(\nabla u,\nabla u)).
\end{equation}
It is easy to see
\begin{equation}\label{non1}
|\mathcal F_1[u]|\le C(|\nabla^2 u|^2+|\nabla u|^4),
\ |\mathcal F_2[u]|\le C(|\nabla^2 u||\nabla u|+|\nabla u|^3),
\end{equation}
where $C>0$ is a constant depending on $\|u\|_{L^\infty(\mathbb R^n)}$.
With the notations as above, (\ref{hfbhm})-(\ref{ic}) can be written as
\begin{equation}
\label{hfbhm1}
({\partial_t}+\Delta^2) u=\mathcal F_1[u]+\nabla\cdot(\mathcal F_2[u])
\ \ {\rm{in}}\  \mathbb R^n\times (0,+\infty); \ u|_{t=0}=u_0.
\end{equation}

The first observation is
\begin{Lemma} For $0<T\le+\infty$,
if $u\in X_T$, then $\mathcal F_1[u]\in Y^1_T, \mathcal F_2[u]\in Y^2_T$.
Moreover,
\begin{equation}\label{map1}
\left\|\mathcal F_1[u]\right\|_{Y^1_T}\le C\left[u\right]_{X_T}^2,
\end{equation}
and
\begin{equation}\label{map2}
\left\|\mathcal F_2[u]\right\|_{Y^2_T}\le C\left[u\right]_{X_T}^2,
\end{equation}
\end{Lemma}
\noindent{\it Proof}. It follows directly from the H\"older inequality. \qed

By the Duhamel formula (2.10), the solution $u$ to (\ref{hfbhm})-(\ref{ic})
 is given by
\begin{equation} \label{duhem}
u=\mathbb Gu_0+\mathbb S(\mathcal F_1[u])+\mathbb S(\nabla\cdot(\mathcal F_2[u])).
\end{equation}

Throughout this section, we denote
$$\hat{u}_0=\mathbb Gu_0.$$
Now we define the mapping operator $\mathbb T$ on $X_{R^4}$ by letting
\begin{equation}\mathbb Tu(x,t)=\hat u_0(x,t)+\mathbb S(\mathcal F_1[u])(x,t)
+\mathbb S(\nabla\cdot(\mathcal F_2[u]))(x,t),\  \ u\in X_{R^4}. \label{t}
\end{equation}

The following property follows directly from Lemma 3.1.
\begin{Lemma} For any $R>0$ and
any initial map $u_0:\mathbb R^n\to N$, $\hat{u}_0\in X_{R^4}$ and
\begin{equation}\label{u0}
\left\|\hat{u}_0\right\|_{L^\infty(\mathbb R^{n+1}_+)}\le C\|u_0\|_{L^\infty(\mathbb R^n)},
\  \ \left[\hat{u}_0\right]_{X_{R^4}}\le C\left[u_0\right]_{\rm{BMO}_R(\mathbb R^n)}.
\end{equation}
\end{Lemma}

For $\epsilon>0$, we define
$$\mathbb B_\epsilon(\hat u_0)
:=\left\{u\in X_{R^4}: \ \|u-\hat u_0\|_{X_{R^4}}\le\epsilon\right\}$$
to be the ball in $X_{R^4}$ with center $\hat u_0$ and radius $\epsilon$.
By the triangle inequality, we have
\beq{}\label{triangle}
\left\|u\right\|_{L^\infty(\mathbb R^{n+1}_+)}\le \epsilon+C\|u_0\|_{L^\infty(\mathbb R^n)},
\ \ \left[u\right]_{X_{R^4}}\le \epsilon+\left[u_0\right]_{{\rm BMO}_R(\mathbb R^n)},
\ \forall u\in \mathbb B_\epsilon(\hat u_0).
\eeq
In particular, we have
\begin{Lemma}For $0<R\le +\infty$,
 if $u_0:\mathbb R^n\to N$ has $[u_0]_{{\rm BMO}_R(\mathbb R^n)}\le\epsilon$, then
\beq{}\label{epsilon}
\left\|u\right\|_{L^\infty(\mathbb R^{n+1}_+)}\le C+\epsilon,
\ \ [u]_{X_{R^4}}\le C\epsilon, \ \forall u\in \mathbb B_\epsilon(\hat u_0)
\eeq
for some $C=C(n, N)>0$.
\end{Lemma}

The proof of Theorem 1.2 is based on the following two lemmas.
\begin{Lemma} There exists $\epsilon_1>0$ 
such that for any $0<R\le+\infty$  if $u_0:\mathbb R^n\to N$ has
$$[u_0]_{{\rm BMO}_R(\mathbb R^n)}\le\epsilon_1,$$
then $\mathbb T$ maps $\mathbb B_{\epsilon_1}(\hat u_0)$ to
$\mathbb B_{\epsilon_1}(\hat u_0)$.
\end{Lemma}
\noindent{\it Proof}.
By (\ref{t}), we have
$$\mathbb T(u)-\hat u_0=\mathbb S(\mathcal F_1[u])
+\mathbb S(\nabla\cdot(\mathcal F_2[u])), \ \
u\in \mathbb B_{\epsilon_1}(\hat u_0).$$
Hence Lemma 4.1, Lemma 4.2, Lemma 5.1, and Lemma 5.2 imply that for any $u\in \mathbb B_{\epsilon_1}(\hat{u}_0)$,
\begin{eqnarray*}
&&\|\mathbb T(u)-\hat u_0\|_{X_{R^4}}\\
&\lesssim& \|\mathbb S(\mathcal F_1[u])\|_{X_{R^4}}+\|\mathbb S(\nabla\cdot(\mathcal F_2[u]))\|_{X_{R^4}}\\
&\lesssim& \|\mathcal F_1[u]\|_{Y^1_{R^4}}+\|\mathcal F_2[u]\|_{Y^2_{R^4}}\\
&\lesssim& \left[u\right]_{X_{R^4}}^2\leq C \epsilon_1^2\le\epsilon_1,
\end{eqnarray*}
provided $\epsilon_1>0$ is chosen to be sufficiently small.
Hence $\mathbb Tu\in \mathbb B_{\epsilon_1}(\hat{u}_0)$.  This completes the proof. \qed

\begin{Lemma} There exist $0<\epsilon_2\le\epsilon_1$ and $\theta_0\in (0,1)$
such that for $0<R\le+\infty$ if $u_0:\mathbb R^n\to N$ satisfies
$$[u_0]_{{\rm BMO}_R(\mathbb R^n)}\le\epsilon_2$$
then $\mathbb T:\mathbb B_{\epsilon_2}(\hat u_0)\to \mathbb B_{\epsilon_2}(\hat u_0)$
is a $\theta_0$-contraction map, i.e.
$$\|\mathbb T(u)-\mathbb T(v)\|_{X_{R^4}}\le\theta_0\|u-v\|_{X_{R^4}},
\ \forall u,v\in \mathbb B_{\epsilon_2}(\hat u_0).
$$
\end{Lemma}
\noindent{\it Proof}. For $u,v\in \mathbb B_{\epsilon_2}(\hat u_0)$, we have
\begin{eqnarray}\label{t_est}
 \left\|\mathbb T u-\mathbb Tv\right\|_{X_{R^4}}
&\leq& \left\|\mathbb S(\mathcal F_1[u]-\mathcal F_1[v])\right\|_{X_{R^4}}
 +\left\|\mathbb S(\nabla\cdot(\mathcal F_2[u]-\mathcal F_2[v]))\right\|_{X_{R^4}}
\nonumber\\
&\lesssim & \left\|\mathcal F_1[u]-\mathcal F_1[v]\right\|_{Y^1_{R^4}}+\left\|\mathcal F_2[u]-\mathcal F_2[v]\right\|_{Y^2_{R^4}}.
\end{eqnarray}
Since
\begin{eqnarray*}
\mathcal F_1[u]-\mathcal F_1[v]&=& -\langle\Delta u, \Delta(D\widetilde\Pi(u))\rangle
+\langle\Delta v, \Delta(D\widetilde\Pi(v))\rangle\\
&=&\langle\Delta (u-v), \Delta(D\widetilde\Pi(u))\rangle
+\langle\Delta v, \Delta(D\widetilde\Pi(u)-D\widetilde\Pi(v))\rangle,
\end{eqnarray*}
we have
\begin{eqnarray*}|\mathcal F_1[u]-\mathcal F_1[v]|&\leq& C [|\Delta(u-v)|(|\Delta u|+|\nabla u|^2+|\Delta v|)\\
&+&|\Delta v|(|\nabla u|+|\nabla v|)|\nabla(u-v)|)]
+ C|\Delta v|(|\nabla^2 u|+|\nabla^2 v|)|u-v|.
\end{eqnarray*}
Hence
\begin{eqnarray}\label{y1_est}
\left\|\mathcal F_1[u]-\mathcal F_1[v]\right\|_{Y^1_{R^4}}
&\leq& C[([u]_{X_{R^4}}+[v]_{X_{R^4}}+[u]_{X_{R^4}}^2)
\|u-v\|_{X_{R^4}} \nonumber\\
&&+[v]_{X_{R^4}}([u]_{X_{R^4}}+[v]_{X_{R^4}})\|u-v\|_{X_{R^4}}]\nonumber\\
&\leq& C\epsilon_2\left\|u-v\right\|_{X_{R^4}},
\end{eqnarray}
where we have used Lemma 5.3 in the last step.

Since
\begin{eqnarray*}
|\mathcal F_2[u]-\mathcal F_2[v]|
&\le &|2(\langle\Delta u, \nabla(D\widetilde\Pi(u))\rangle-
\langle\Delta v, \nabla(D\widetilde \Pi(v))\rangle)|
\\
&+&|\nabla(D^2\widetilde\Pi(u)(\nabla u,\nabla u)-D^2\widetilde\Pi(v)(\nabla v,\nabla v))|\\
&\leq & C [|\nabla u||\Delta (u-v)|+|\Delta v|(|u-v|+|\nabla(u-v))|]\\
&+&C[|\nabla u|(|\nabla u|+|\nabla v|)|\nabla(u-v)|+(|\nabla^2 u|+|\nabla^2 v|)|\nabla (u-v)|]\\
&+&C[(|\nabla u|+|\nabla v|)|\nabla^2(u-v)|+|\nabla v|^2|\nabla(u-v)|+|\nabla v||\nabla^2 v||u-v|],
\end{eqnarray*}
we have
\begin{eqnarray}\label{y2_est}
\left\|\mathcal F_2[u]-\mathcal F_2[v]\right\|_{Y^2_{R^4}}
&\leq& C([u]_{X_{R^4}}+[v]_{X_{R^4}}+[u]_{X_{R^4}}^2+[v]_{X_{R^4}}^2)
\|u-v\|_{X_{R^4}}\nonumber\\
&\le& C\epsilon_2\|u-v\|_{X_{R^4}}.
\end{eqnarray}
Putting (\ref{y1_est}) and (\ref{y2_est}) into (\ref{t_est}), we obtain
$$\|\mathbb Tu-\mathbb Tv\|_{X_{R^4}}
\leq C\epsilon_2\|u-v\|_{X_{R^4}}\leq\theta_0 \|u-v\|_{X_{R^4}}$$
for some $\theta_0=\theta_0(\epsilon_2)\in (0,1)$, provided $\epsilon_2>0$ is chosen to be  sufficiently small.
This completes the proof.
\qed\\

\noindent{\bf Proof of Theorem 1.2}. It follows from Lemma 5.4 and Lemma 5.5, and the fixed point theorem
that there  exists $\epsilon_0>0$ such that
for $0<R\le+\infty$ if $[u_0]_{{\rm BMO}_R(\mathbb R^n)}\le\epsilon_0$,
then there exists a unique $u\in X_{R^4}$ such that
$$u=\hat{u}_0+\mathbb S(\mathcal F[u])\ \ {\rm{on}}\ \mathbb R^n\times[0,R^4),$$
or equivalently 
$$u_t+\Delta^2 u = \mathcal F[u] \ {\rm{on}}\ \mathbb R^n\times (0,R^4);
\ u\big|_{t=0}=u_0.$$ 

Now we need to show $u(\mathbb R^n\times [0,R^4])\subset N$. First, observe that
Lemma 2.1 implies that for any $x\in\mathbb R^n$ and $t\le \frac{R^4}{K^4}$, 
\begin{eqnarray*}
{\rm{dist}}(u(x,t),N)&\leq &{\rm{dist}}(\hat{u}_0(x,t),N)
+\|u-\hat{u}_0\|_{L^\infty(\mathbb R^n\times [0,R^4])}\\
&\leq& \delta +K^n\left[u_0\right]_{{\rm BMO}_R(\mathbb R^n)}+\epsilon_0\\
&\leq &\delta+(1+K^n)\epsilon_0\le\delta_N,
\end{eqnarray*}
provide $\delta\le \frac{\delta_N}2$ and $\epsilon_0\le \frac{\delta_N}{2(1+K^n)}$.  
This yields $u(\mathbb R^n\times [0,
\frac{R^4}{K^4}])\subset N_{\delta_N}$, the $\delta_N$-neighborhood of $N$. 
This and the definition of $\widetilde\Pi(\cdot)$ imply that  $
\widetilde\Pi(u)\equiv\Pi(u)$.

Set $Q(y)=y-\Pi(y)$ for $y\in N_{\delta_N}$,
and $\rho(u)=\frac12|Q(u)|^2$. Then direct calculations imply that for any $y\in N_{\delta_N}$,
$$D Q(y)(v)=({\rm{Id}}-D\Pi(y))(v), \ \forall v\in \mathbb R^l,$$
and
$$D^2 Q(y)(v,w)=-D^2\Pi(y)(v,w), \ \forall v, w\in\mathbb R^l.$$
Observe that $\mathcal F[u]$ can be rewritten as
\begin{eqnarray*}
&&\mathcal F[u]\\
&=&\Delta(D^2\Pi(u)(\nabla u,\nabla u))+\nabla\cdot(D^2\Pi(u)(\Delta u,\nabla u)) +
D^2\Pi(u)(\nabla\Delta u,\nabla u).
\end{eqnarray*}
Direct calculations imply
\begin{eqnarray}\label{rho_eqn}
&&(\partial_t+\Delta^2)Q(u)\nonumber\\
&=&D Q(u)(\partial_t u+\Delta^2 u)\nonumber\\
&-&\left[D^2\Pi(u)(\nabla\Delta u, \nabla u)+\nabla\cdot(D^2\Pi(u)(\Delta u, \nabla u))
+\Delta(D^2\Pi(u)(\nabla u,\nabla u))\right]\nonumber\\
&=&D Q(u)(\mathcal F[u])-\mathcal F[u]\nonumber\\
&=&-D\Pi(u)(\mathcal F[u]).
\end{eqnarray}
Multiplying both sides of (\ref{rho_eqn}) by $Q(u)$ and integrating over $\mathbb R^n$, we
obtain
\begin{eqnarray}\label{gronwall}
\frac{d}{dt} \int_{\mathbb R^n}\rho(u)+\frac12\int_{\mathbb R^n}|\Delta (Q(u))|^2
&=&-\frac12\int_{\mathbb R^n}\langle D\Pi(u)(\mathcal F[u]), Q(u)\rangle\nonumber\\
&=&0,
\end{eqnarray}
where we have used the fact that $Q(u)\perp T_{\Pi(u)}N$ and 
$D\Pi(u)(\mathcal F[u])\in T_{\Pi(u)}N$ in the last step. 

Since $\rho(u)|_{t=0}=0$, integrating (\ref{gronwall}) from $0$ to $\frac{R^4}{K^4}$ implies $\rho(u)\equiv 0$
on $\mathbb R^n\times [0,\frac{R^4}{K^4}]$. 
 Thus $u(\mathbb R^n\times [0,\frac{R^4}{K^4}])\subset N$.
Repeating the same argument for $t\in [\frac{R^4}{K^4}, R^4]$ yields
$u(\mathbb R^n\times [\frac{R^4}{K^4}, R^4])\subset N$.
This completes the proof of Theorem 1.2. \qed\\

\noindent{\bf Proof of Theorem 1.3}. It follows directly from Theorem 1.2 with
$R=+\infty$. \qed

\setcounter{section}{5} \setcounter{equation}{0}
\section{Proof of Theorem 1.4 and 1.5}

This section is devoted to the proof of both Theorem 1.4 and 1.5. Since the argument is similar
to that of Theorem 1.2, we will only sketch it here. 

Let $\mathcal H[u]$ denote the right hand side of (\ref{int_bih}).
Then we have
$$\mathcal H[u]=\mathcal F_1[u]+\nabla\cdot\mathcal F_2[u]+\mathcal F_3[u],$$
where $\mathcal F_1[u]$ and $\mathcal F_2[u]$ are  given by (5.2), while
\begin{eqnarray}
\mathcal F_3[u]&=&D\widetilde\Pi(u)[D^2\widetilde\Pi(u)(\nabla u,\nabla u)\cdot D^3\widetilde\Pi(u)(\nabla u,\nabla u)]
\nonumber\\
&&+2D^2\widetilde\Pi(u)(\nabla u,\nabla u)\cdot D^2\widetilde\Pi(u)(\nabla u, \nabla(D\widetilde\Pi(u))).
\end{eqnarray}

It is clear that $u\in X_{R^4}$ solves (\ref{int_bih})-(\ref{ic1}) iff 
\beq{} u=\mathbb G u_0+\mathbb S(\mathcal F_1[u])+\mathbb S(\nabla\cdot\mathcal F_2[u])+\mathbb S(\mathcal F_3[u]).
\eeq
Since $\mathcal F_3[u]$ satisfies
\beq{}
|\mathcal F_3[u]|\le C|\nabla u|^4,
\eeq
for some $C>0$ depending on $\|u\|_{L^\infty(\mathbb R^n)}$, it is easy to check\\
\noindent{\bf Claim 1}. {\it For $0<R\le +\infty$, if $u\in X_{R^4}$, then $\mathcal F_3[u]\in Y^1_{R^4}$ and
\beq{}
\|\mathcal F_3[u]\|_{Y^1_{R^4}}\le C\left[u\right]_{X_{R^4}}^4.
\eeq}
This claim and Lemma Lemma 4.1 then imply\\
\noindent{\bf Claim 2}. {\it For $0<R\le +\infty$, if $u\in X_{R^4}$, then $\mathbb S(\mathcal F_3[u])\in X_{R^4}$ and
\beq{}
\|\mathbb S(\mathcal F_3[u])\|_{X_{R^4}}\le C\left[u\right]_{X_{R^4}}^4.
\eeq}
Now if define the mapping operator $\widetilde{\mathbb T}$ on $X_{R^4}$ by
\beq{}
\widetilde{\mathbb T}[u]:=
\mathbb G u_0+\mathbb S(\mathcal F_1[u])+\mathbb S(\nabla\cdot\mathcal F_2[u])+\mathbb S(\mathcal F_3[u]),
\eeq
then Claim 1, Claim 2, and Lemma 5.4 imply\\
\noindent{\bf Claim 3}. {\it There exists $\epsilon_3>0$ such that for $0<R\le+\infty$,
if $u_0:\mathbb R^n\to N$ has $[u_0]_{{\rm BMO}_R(\mathbb R^n)}\le\epsilon_3$, then
$\widetilde{\mathbb T}$ maps $\mathbb B_{\epsilon_3}(\hat{u}_0)$ to $\mathbb B_{\epsilon_3}(\hat{u}_0)$.}

We need to show $\widetilde{\mathbb T}:\mathbb B_{\epsilon_3}(\hat{u}_0)\to \mathbb B_{\epsilon_3}(\hat{u}_0)$
is a contraction map. To see this, observe that direct calculations imply
that for any $u, v\in \mathbb B_{\epsilon_3}(\hat{u}_0)$,
\begin{eqnarray}
&&|\mathcal F_3[u]-\mathcal F_3[v]|\\
&\leq& C[|u-v||\nabla u|^4+|\nabla(u-v)|(|\nabla v|^3+|\nabla v||\nabla u|^2+|\nabla v|^2|\nabla u|+|\nabla u|^3)]
\nonumber
\end{eqnarray}
for some $C>0$ depending only $\max\{\|u\|_{L^\infty(\mathbb R^n)}, \|v\|_{L^\infty(\mathbb R^n)}\}$.
Hence, combined with the proof of Lemma 5.5,  we obtain\\
\noindent{\bf Claim 4}. {\it There exists $\epsilon_3>0$ such that for $0<R\le+\infty$,
if $u_0:\mathbb R^n\to N$ has $[u_0]_{{\rm BMO}_R(\mathbb R^n)}\le\epsilon_3$, then
\beq{}
\|\widetilde{\mathbb T}[u]-\widetilde{\mathbb T}[v]\|_{X_{R^4}}
\le C\epsilon_3 \|u-v\|_{X_{R^4}}, \ \forall u, v\in \mathbb B_{\epsilon_3}(\hat{u}_0).
\eeq
}

Now we can complete the proof of Theorem 1.4 as follows.

\noindent{\bf Completion of proof of Theorem 1.4}:  Similar to Theorem 1.2, it follows from Claim 3 and Claim 4
and the fixed point theorem that  there  exists $\epsilon_0>0$ such that
for $0<R\le+\infty$ if $[u_0]_{{\rm BMO}_R(\mathbb R^n)}\le\epsilon_0$,
then there exists a unique $u\in X_{R^4}$ that solves (\ref{int_bih})-(\ref{ic1}):
$$u_t+\Delta^2 u = \mathcal H[u] \ {\rm{on}}\ \mathbb R^n\times (0,R^4);
\ u\big|_{t=0}=u_0.$$ 
The same argument as in Theorem 1.2 implies $u(\mathbb R^n\times [0,\frac{R^4}{K^4}])\subset N_{\delta_N}$.
Hence $\widetilde{\Pi}(u)\equiv\Pi(u)$ on $\mathbb R^n\times [0,\frac{R^4}{K^4}]$.
Moreover, the same calculation as in (5.15) implies
\beq{}
(\partial_t+\Delta^2)(u-D\Pi(u))=-D\Pi(u)(\mathcal H[u]),
\eeq
and it follows that for $0\le t\leq\frac{R^4}{K^4}$, 
\beq{}
\frac{d}{dt}\int_{\mathbb R^n}|u-D\Pi(u)|^2+\int_{\mathbb R^n}|\Delta(u-D\Pi(u))|^2=0.
\eeq
This, combined with $|u-D\Pi(u)|^2\big|_{t=0}=0$, implies that $u(\mathbb R^n\times [0,\frac{R^4}{K^4}])\subset N$.
Repeating the same argument then implies $u(\mathbb R^n\times [0,R^4])\subset N$. 
The proof is complete. \qed\\

\noindent{\bf Proof of Theorem 1.5}. It follows directly from Theorem 1.4 with $R=+\infty$. \qed\\

\noindent{\bf Acknowledgements}. The work is partially supported by NSF grant 0601182.
This work is carried out while the author is visiting IMA as a New Directions Research Professorship. 
The author wishes to thank IMA for providing both the financial support
and the excellent research environment.


\begin{thebibliography}{99}

\bibitem{A} G. Angelesberg, {\em A monotonicity formula for stationary biharmonic maps.}
Math. Z. 252 (2006), no. 2, 287-293.

\bibitem{CWY} A. Chang, L. Wang, P. Yang, {\em A regularity theory of biharmonic maps.}
Comm. Pure Appl. Math. 52:1113-1137.

\bibitem{G} A. Gastel, {\em The extrinsic polyharmonic map heat flow in the critical dimension.} Adv. Geom. 6 (2006), no. 4, 501-521.

\bibitem{HW} M. Hong, C. Wang, {\em Regularity and relaxed problems of minimizing biharmonic maps into spheres.}  Calc. Var. Partial Differential Equations 23 (2005), no. 4, 425-450. 

\bibitem{HW1} T. Huang, C. Wang, {\em Well-posedness for the heat flow of polyharmonic maps with rough initial data}.
In preparation.

\bibitem{KL} H. Koch, T. Lamm, {\em Geometric flows with rough initial data}. arXiv:
0902.1488v1, 2009.

\bibitem{KT} H. Koch, D. Tataru, {\em Well-posedness for the Navier-Stokes equations.}
Adv. Math. 157 (2001), no. 1, 22-35.

\bibitem{L1} T. Lamm, {\em Biharmonischer W\"armefluss.} Diplomarbeit Universit\"at Freiburg (2001).

\bibitem{L2} T. Lamm, {\em Heat flow for extrinsic biharmonic maps with small initial energy.} Ann. Global Anal. Geom. 26 (2004), no. 4, 369-384.

\bibitem{L3} T. Lamm, {\em Biharmonic map heat flow into manifolds of nonpositive curvature.} Calc. Var. Partial Differential Equations 22 (2005), no. 4, 421-445.

\bibitem{LR} T. Lamm, T. Riviere, {\em Conservation laws for fourth order systems in four dimensions.} Comm. Partial Differential Equations 33 (2008), no. 1-3, 245-262. 

\bibitem{Sch} C. Scheven, {\em Dimension reduction for the singular set of biharmonic maps.} Adv. Calc. Var. 1 (2008), no. 1, 53-91.

\bibitem{Stein} E. Stein, "Harmonic Analysis", Princeton Mathematical Series, Vol. 43, Princeton
University Press, Princeton, 1993.

\bibitem{Struwe} M. Struwe, {\em Partial regularity for biharmonic maps, revisited.}
Calc. Var. Partial Differential Equations 33 (2008), no. 2, 249-262.

\bibitem{Str} P. Strzelecki, {\em On biharmonic maps and their generalizations. }
Calc. Var. Partial Differ. Equ. 18(4), 401-432 (2003).

\bibitem{W1} C. Wang, {\em Biharmonic maps from $\Bbb{R}^{4}$ into a Riemannian manifold.} Math. Z. 247:65-87.

\bibitem{W2} C. Wang, {\em Stationary biharmonic maps from $\Bbb{R}^{m}$ into a Riemannian manifold.} Comm. Pure Appl. Math. 57:419-444.

\bibitem{W3} C. Wang, {\em Remarks on biharmonic maps into spheres.}
 Calc. Var. Partial Differential Equations 21:221-242.

\bibitem{W4} C. Wang, {\em Heat flow of biharmonic maps in dimensions four and its application.}  Pure Appl. Math. Q. 3 (2007), no. 2, part 1, 595-613.

\bibitem{W5} C. Wang, {\em Well-posedness for the heat flow of harmonic maps and the liquid
crystal flow with rough initial data }. Preprint (2010).

\end{thebibliography}
\end{document}